\documentclass[11pt]{article}

\usepackage{latexsym,color,amsmath,amsthm,amssymb,amscd,amsfonts}

\usepackage{graphicx}
\usepackage{cancel}
\usepackage{stmaryrd}
\usepackage{commath}
\usepackage{mathtools}
\usepackage{amssymb}
\usepackage{enumerate}
\usepackage{enumitem}
\usepackage{bbm}
\usepackage{hyperref,mathtools}
\usepackage[title]{appendix}

\newtheorem{conj}{Conjecture}[section]
\newtheorem{theorem}[conj]{Theorem}

\newtheorem{remark}[conj]{Remark}

\newtheorem{ex}[conj]{Example}

\newcommand\independent{\protect\mathpalette{\protect\independent}{\perp}} 
\def\independent#1#2{\mathrel{\rlap{$#1#2$}\mkern2mu{#1#2}}}

\newcommand\restr[2]{{
  \left.\kern-\nulldelimiterspace 
  #1 
  \littletaller 
  \right|_{#2} 
  }}

\newcommand{\littletaller}{\mathchoice{\vphantom{\big|}}{}{}{}}

\date{}

\author{Heshan Aravinda}

\title{A note on the maximum probability of ultra log-concave distributions}

\begin{document}

\maketitle
\begin{abstract} \noindent
Jakimiuk et al. (2024) have proved that, if $X$ is an ultra log-concave random variable with integral mean, then
$$\max_n \mathbb{P}\{X=n\} \geq \max_n \mathbb{P} \{Z=n\}\,,$$ where $Z$ is a Poisson random variable with the parameter $\mathbb{E}[X]$. In this note, we show that this inequality does not always hold true when $X$ is ultra log-concave with $\mathbb{E}[X]>1$.
\end{abstract}
\section{Introduction}
\label{intro}
A discrete random variable $X$ with the probability mass function $p$ is said to be \textit{log-concave} if the support of $X$ forms an interval of integers and 
$$p^2(n) \geq p(n-1)\,p(n+1)\,\,\,\,\,\,\,\text{for all $n \in \mathbb{Z}$.}$$
A stronger notion, analogous to strong log-concavity in the case of continuous
 random variables, is that of \textit{ultra-log-concavity}. This notion of ultra log-concavity arises from the search for a theory of negative dependence. As argued by Pemantle, a theory of negative dependence has long been desired in probability and statistical physics, in analogy with the theory of positive dependence (see \cite{Pem}). Recall a random variable $X$ taking values in $\{0,1,2,3,...\}$ is defined to be ultra log-concave if its probability mass function $p$ with respect to the law of a Poisson distribution, is log-concave, i.e. $$p\,^2(n) \geq (1+ \frac{1}{n})\,p(n+1)\,p(n-1)\,\,\,\,\,\,\text{for all $n = 1,2,3,\dots$}$$ Poisson distribution as well as all Bernoulli sums are ultra-log-concave. A non-trivial property of the class of ultra-log-concave distributions, is that it is closed under convolution \cite{L}. Some recent works concerning this class of random variables include ultra log-concavity of discrete order statistics \cite{BDP}, maximum entropy property of Poisson \cite{JKM}, concentration inequalities \cite{HesMarMel} and relative log-concavity ordering \cite{WW}. \vskip2mm
Our main interest is the following inequality established by Jakimiuk et al. (see \cite[Theorem 2, Corollary 3]{JMNS}). See also \cite{H} for a related result by Hoeffding for Bernoulli sums.
\begin{theorem}
    \label{JMN_ULC}
    If $X$ is an ultra log-concave random variable with integral mean, then
    \begin{equation}     \label{JMN}
    \max_n \mathbb{P}\{X=n\} \geq \max_n \mathbb{P} \{Z=n\}\,,
     \end{equation}
     where $Z \sim \text{Pois}(\mathbb{E}[X])$. 
\end{theorem}
Therein, the authors mentioned that it would be natural to expect the inequality (\ref{JMN}) to hold under relaxed assumptions than $\mathbb{E}[X]$ being an integer, particularly when $X$ is ultra log-concave with $\mathbb{E}[X]\geq 1$. In this note, we show that this is not possible. \vskip2mm

\begin{remark} Note that for certain naturally occurring ULC distributions, one can expect a strict inequality to hold.
\begin{enumerate}
\item \textbf{Zero-truncated Poisson distribution (ZTP).}  An integer valued random variable $X$ is said to be ZTP with the parameter $\lambda>0$, if its probability mass function is of the following form:
$$\mathbb{P}\{X=n\} = \dfrac{\lambda^n}{(e^\lambda-1)\,n!}\,\,,\,\,\,\,n=1,2,3,\dots$$
Clearly, $X$ is ULC. Moreover, $\mathbb{E}[X] = \dfrac{\lambda\,e^\lambda}{e^\lambda -1}>1$ for all $\lambda>0$. Let $\lambda^\prime = \lambda \frac{e^\lambda}{e^\lambda-1}$. Consider $Z \sim \text{Pois} (\mathbb{E}[X])$. Let us verify the inequality: $\displaystyle \max_n \mathbb{P}\{X=n\} > \displaystyle \max_n \mathbb{P}\{Z=n\}$.\vskip2mm
Notice that the ratio $\frac{\mathbb{P}\{X=n+1\}}{\mathbb{P}\{X=n\}}= \frac{\lambda}{n+1} \geq 1$ when $n \leq \lambda -1$, and $\frac{\mathbb{P}\{X=n+1\}}{\mathbb{P}\{X=n\}}<1$ when $n > \lambda -1$, so the mode of $X$ is $\lfloor \lambda \rfloor$ except when $0<\lambda<1$ in which case the mode is $n=1$. For $Z$ Poisson, the mode is $\lfloor \lambda^\prime \rfloor$. Therefore, we show that $\mathbb{P}\{X= \lfloor \lambda \rfloor\} > \mathbb{P}\{Z= \lfloor \lambda^\prime \rfloor\}$ for $\lambda>0$. The case $0<\lambda<1$ is straightforward. Suppose $\lambda \geq 1$. Then
\begin{align*}
\mathbb{P}\{X= \lfloor \lambda \rfloor\} > \mathbb{P}\{Z= \lfloor \lambda^\prime \rfloor\} & \iff  \dfrac{\lambda ^{\lfloor \lambda \rfloor}}{(e^\lambda -1)\, \lfloor \lambda \rfloor!}> \dfrac{{\lambda^\prime} ^{\lfloor \lambda^\prime \rfloor}}{e^{\lambda^\prime}\, \lfloor \lambda^\prime \rfloor!}\\
&\iff  1> \dfrac{(e^\lambda -1)\,\lfloor \lambda \rfloor!}{ \lambda ^{\lfloor \lambda \rfloor}} \cdot \dfrac{{\lambda^\prime} ^{\lfloor \lambda^\prime \rfloor}}{e^{\lambda^\prime}\, \lfloor \lambda^\prime \rfloor!} \\
\end{align*}
To this end, define $f(x) = \frac{e^x\,\lfloor x \rfloor !}{x^{\lfloor x \rfloor}}$ on $[1, \infty)$. Since $\lambda < \lambda^\prime$, it suffices to show that $f$ is increasing. Observe that for any integer $n\geq 1$, $f(x) = \frac{e^x n!}{x^n}$ on the interval $(n, n+1)$. Therefore, $f^\prime (x)= \frac{e^x n!}{x^{n+1}} (x-n)$, which is positive on $(n, n+1)$. Since $f(n)= \lim_{x \to n^{+}} f(x)  = \frac{e^n n!}{n^n} = \frac{e^n (n-1)!}{n^{n-1}}= \lim_{x \to n^{-}} f(x)$, the function $f$ is also continuous at $x=n$ which ensures that there are no jumps at integer points.

\item \textbf{Binomial distribution.}  Let $X \sim \text{Bin}(m,p)$ with $\mathbb{E}[X]\geq 1$. Fix an integer $k\geq 1$. For any integer $m \geq k+1$, choose $p \in [\frac{k}{m}, \frac{k+1}{m})$ so that $\lfloor mp \rfloor = k$. It is well-known that the mode of $X$ is $\lfloor (m+1)p \rfloor$. Therefore, we show that $\mathbb{P} \{X= \lfloor (m+1)p \rfloor\} > \mathbb{P} \{Z= \lfloor mp \rfloor\}$, where $\lfloor mp \rfloor$ is the mode of $Z$. \vskip2mm 

There are two cases: $\lfloor (m+1)p \rfloor=k$ for $\frac{k}{m}\leq p < \frac{k+1}{m+1}$ and  $\lfloor (m+1)p \rfloor = k+1$ for $\frac{k+1}{m+1}\leq p < \frac{k+1}{m}$. Let us first treat the case $\lfloor (m+1)p \rfloor=k$. For this case, after simplifying, the desired inequality becomes

$$ \dfrac{m!}{(m-k)!}\, (1-p)^{m-k} > e^{-mp} \, m^k.$$

In particular, we wish to show that for any integer $m\geq k+1$ and $p$ satisfying $\frac{k}{m}\leq p < \frac{k+1}{m+1}$
$$ \dfrac{m!}{(m-k)!}\, (1-p)^{m-k} > e^{-mp} \, m^k.$$

Proceed by induction on $k$. When $k=1$, the above inequality gets reduced to
$$(1-p)^{m-1} > e^{-mp}.$$
Equivalently, after taking the natural logarithms on both sides, we show that for any integer $ m\geq 2$ and $\frac{1}{m} \leq p < \frac{2}{m+1}$
$$f_m(p) := (m-1)\,\ln{(1-p)} +mp >0.$$
Differentiating $f_m$ w.r.t $p$ gives $f_m^\prime (p) = \frac{1-mp}{1-p}$, which is non-positive on $[\frac{1}{m}, \frac{2}{m+1})$. Therefore, $f_m(p) \geq f_m(\frac{2}{m+1}) = \frac{2m}{1+m} + (m-1) \ln{\left(\frac{m-1}{m+1}\right)}>0$. We claim that $f_m(\frac{2}{m+1})>0$ for $m \geq 2$. Indeed, after taking the derivative of $f_m(\frac{2}{m+1})$ w.r.t $m$, we get $f_m^\prime(\frac{2}{m+1}) = \frac{2}{m+1}+ \frac{2}{(m+1)^2} +\ln{\left(\frac{m-1}{m+1}\right)}$. By letting $x = \frac{2}{m+1}$ in the Maclaurin series $\ln{(1-x)} = -x - \frac{x^2}{2} + O(x^3)$ on $[-1, 1)$, we have $\ln{\left(\frac{m-1}{m+1}\right)}< -\frac{2}{m+1} - \frac{2}{(m+1)^2}$, which shows that $f_m(\frac{2}{m+1})$ is decreasing for $m \geq 2$. Finally, we conclude since $f_m(\frac{2}{m+1}) \to 0$ as $m \to \infty$.
\vskip2mm
Let us now assume the inequality holds for $m\geq k+1$ and $\frac{k}{m}\leq p <\frac{k+1}{m+1}$. Then, for $m\geq k+2$ and $\frac{k+1}{m}\leq p < \frac{k+2}{m+1}$
\begin{align*}
\dfrac{m!}{(m-k-1)!}\, (1-p)^{m-k-1} - e^{-mp} \, m^{k+1} & = \frac{m-k}{1-p}\dfrac{m!}{(m-k)!} (1-p)^{m-k} - e^{-mp} m^{k+1}\\
&> \dfrac{m-k}{1-p}\,e^{-mp} \, m^k - e^{-mp} \, m^{k+1}\\
& = e^{-mp} \, m^k \left(\dfrac{m-k}{1-p} - m\right)>0.
\end{align*} 
Here, the first inequality is due to the induction hypothesis and the second inequality is due to the fact that $mp\geq k+1$.

The remaining case follows from the proof of the first case, since $\lfloor (m+1)p \rfloor = k+1 \iff \frac{k+1}{m+1}\leq p < \frac{k+1}{m}$ for $m\geq k+1$, and
$$ \max_n \mathbb{P}\{X=n\} =  \mathbb{P}\{X= \lfloor (m+1)p \rfloor \coloneqq k+1 \}\geq \mathbb{P}\{X= k\}> \mathbb{P}\{Z= \lfloor mp \rfloor\coloneqq k \} .$$

\end{enumerate}

\end{remark}

Let us also state an entropic version of the inequality (\ref{JMN}): $H_{\infty}(X) \leq H_{\infty}(Z)$ for $X$ - ultra log-concave with integral mean. Here,  $ H_{\alpha}(X) =\frac{-1}{\alpha-1}\log \sum_{n \in \mathbb{Z}}\,p(n)^{\alpha},\, \alpha \neq 1$ is defined as the \textit{the R\'{e}nyi entropy of order $\alpha$}, with $H_{\infty}(X)$ being identified as the limiting case $\alpha \to \infty$; $H_\infty (X) = - \log \max_n p(n)$, referred to as the \textit{min-entropy}. As mentioned in \cite{JMNS}, it follows from \cite[Theorem 1.13]{MM2} that for any ultra log-concave random variable $X$, $H_{\alpha}(X) \leq H_{\alpha}(Z)$ for $\alpha \leq 1$. However, one cannot extend this for $\alpha>1$ (see \cite[Proposition 5.3]{MM2}).

\section{An example: ultra log-affine distributions}
This section is devoted to the construction of a distribution for which the natural extension of (\ref{JMN}) is not valid. The idea is to consider a compactly supported ultra log-affine random variable with the parameter $\lambda$, i.e. an integer-valued random variable $X$ whose mass function defined as $p(n) = C\, \frac{\lambda^n}{n!}\,,\,n \in \{K,K+1,\dots,N\}$, where $\,0\leq K \leq N$ and $C$ is the normalizing constant. Such random variables have been identified as an extremal case in a recent work by \cite{MarMel}, where the authors developed a localization-type technique for discrete log-concave measures.
\begin{ex}
Fix $N\geq 5$, and let $X$ be ultra log-affine with the parameter $\frac{3}{2}\leq\lambda \leq 2$, supported on $\{0,1,2,\dots,N\}$. Then
\begin{equation}
\label{counter}
    \max_n \mathbb{P}\{X=n\} < \max_n \mathbb{P} \{Z=n\}\,,
\end{equation}
where $Z \sim \text{Pois} (\mathbb{E}[X]) $
\end{ex}
\begin{proof}
Since $X$ is ultra log-affine, its probability mass function $p$ takes the form $p(n) = C\, \frac{\lambda^n}{n!}\,,\,n \in \{0,1,2,\dots,N\}$, where $C = \dfrac{1}{\sum_{n=0}^N \frac{\lambda^n}{n!}}$. Then,
    \begin{align*}
        \mathbb{E}[X] = C\sum_{n=0}^N n\,\dfrac{\lambda^n}{n!} = \dfrac{\lambda\,\displaystyle \sum_{n=0}^{N-1} \dfrac{\lambda^n}{n!}}{\displaystyle \sum_{n=0}^N \dfrac{\lambda^n}{n!}}.
    \end{align*}
Note that $1<\mathbb{E}[X]<2$. Indeed since
    \begin{align*}
       \mathbb{E}[X]>1 \iff \lambda \,\sum_{n=0}^{N-1} \dfrac{\lambda^n}{n!} - \sum_{n=0}^{N} \dfrac{\lambda^n}{n!} >0 & \iff \sum_{n=0}^{N-1} (\lambda -1)\,\dfrac{\lambda^n}{n!} - \dfrac{\lambda^N}{N!}>0\\
        & \iff  \sum_{n=0}^{N-2}(\lambda -1)\,\dfrac{\lambda^n}{n!} + (\lambda -1)\,\dfrac{\lambda^{N-1}}{(N-1)!}  - \dfrac{\lambda^N}{N!} >0\\
        & \iff \sum_{n=0}^{N-2}(\lambda -1)\,\dfrac{\lambda^n}{n!} + (\lambda -1)\,\dfrac{\lambda^{N-1}}{N!} \left(N-\frac{\lambda}{\lambda -1}\right) >0.
    \end{align*}
    The last inequality is true since $N-\frac{\lambda}{\lambda -1} = N-1 - \frac{1}{\lambda-1} \geq N-3>0$ for $3/2\leq \lambda \leq 2$. And, it is easy to see $\mathbb{E}[X]<2$. \vskip1mm
    Let us now compute $\displaystyle \max_n \mathbb{P}\{X=n\}$. Since $\dfrac{\mathbb{P}\{X=n\}}{\mathbb{P}\{X=n-1\}} = \dfrac{\lambda}{n}$ for $1\leq n\leq N$, we have 
    $$\max_n \mathbb{P}\{X=n\} =  \mathbb{P}\{X=1\} = \dfrac{\lambda}{\displaystyle \sum_{n=0}^{N} \dfrac{\lambda^n}
        {n!} }.$$
        On the other hand, if $Z \sim \text{Pois} (\mathbb{E}[X])$, then $\displaystyle \max_n \mathbb{P} \{Z=n\}$  occurs when $n =1$. To see why, consider the ratio
        $$\dfrac{\mathbb{P}\{Z=n\}}{\mathbb{P}\{Z=n-1\}} = \dfrac{\mathbb{E}[X]}{n}\,\,\,\,\,\,\,\,\text{for}\,\,\,1\leq n \leq N,$$
        which is less than $1$ if $n\geq 2$ and greater than $1$ if $n=1$. Therefore,
        $$\max_n \mathbb{P}\{Z=n\} =  \mathbb{P}\{Z=1\} = e^{-\mathbb{E}[X]}\,\mathbb{E}[X].$$
        We wish to show
        \begin{equation}
    \max_n \mathbb{P}\{X=n\}<\max_n \mathbb{P}\{Z=n\} \iff \dfrac{\lambda}{\displaystyle \sum_{n=0}^{N} \dfrac{\lambda^n}
        {n!} } < e^{-\mathbb{E}[X]}\,\mathbb{E}[X].
        \end{equation}
        After substituting for $\mathbb{E}[X]$ and rearranging the terms, inequality (3) is equivalent to
        \begin{equation} \label{two_point}
            \exp \left( \dfrac{-\lambda\, \displaystyle \sum_{n=0}^{N-1} \dfrac{\lambda^n}
        {n!} }{\displaystyle \sum_{n=0}^{N} \dfrac{\lambda^n}
        {n!}}\right)\, \sum_{n=0}^{N-1} \dfrac{\lambda^n}
        {n!}>1 \,\,\,\,\,\, \text{for all $2\geq \lambda \geq \frac{3}{2}$ and $N\geq 5$.}
        \end{equation}
        \vskip4mm \noindent
        \subsection{A Two-point Inequality}  It remains to show inequality (\ref{two_point}). First, we need the following:\vskip3mm \noindent
\textbf{\underline{Fact I:}}  $ \displaystyle \sum\limits_{n=0}^{N-1} \dfrac{\lambda^n}{n!} \geq e^\lambda\Big(1 -\frac{\lambda^N}{e^\lambda\,N!}\ \frac{N+1}{N+1-\lambda}\Big).$
\begin{proof}
If $ n\geq N$ then $\ n! = N!(N+1)(N+2)\dots n \geq N!(N+1)^{n-N}$, leading to
$$
  \sum_{n=N}^\infty \frac{\lambda^n}{n!} \leq
   \sum\limits_{n=N}^\infty \frac{\lambda^n}{N!\,(N+1)^{n-N}} =
\frac{\lambda^N}{N!}\sum\limits_{n=N}^\infty \Big(\frac{\lambda}{N+1}\Big)^{\large \!n-N} \!\!=
\frac{\lambda^N}{N!}\ \frac{N+1}{N+1-\lambda},
$$
where the final sum is evaluated as a geometric sequence. This gives us:
$$
  \sum\limits_{n=0}^{N-1} \frac{\lambda^n}{n!} = e^\lambda - \sum\limits_{n=N}^\infty \frac{\lambda^n}{n!} 
 \geq e^\lambda -\frac{\lambda^N}{N!}\ \frac{N+1}{N+1-\lambda} =
e^\lambda\Big(1 -\frac{\lambda^N}{e^\lambda\,N!}\ \frac{N+1}{N+1-\lambda}\Big) .
$$
\end{proof}
\noindent
\textbf{\underline{Fact II:}} $\exp \left( \dfrac{-\lambda\, \displaystyle \sum_{n=0}^{N-1} \dfrac{\lambda^n}{n!} }{\displaystyle \sum_{n=0}^{N} \dfrac{\lambda^n}{n!}}\right) >  e^{-\lambda}\ \left(1+\lambda\ {\frac{\lambda^N}{e^\lambda \,N!}}\right).$
\begin{proof}
\begin{align*}
\exp \left( \frac{-\lambda\, \displaystyle \sum_{n=0}^{N-1} \frac{\lambda^n}
        {n!} }{\displaystyle \sum_{n=0}^{N} \frac{\lambda^n}
        {n!}}\right)
= \exp \left(\dfrac{-\lambda \displaystyle \sum_{n=0}^{N} \frac{\lambda^n}{n!} + \lambda\ \frac{\lambda^N}{ N!}}{\displaystyle \sum_{n=0}^N \frac{\lambda^n}{n!}}\right) &= e^{-\lambda}\ \exp \left(\dfrac{ \displaystyle \lambda\  \frac{\lambda^N}{ N!}}{\displaystyle \sum_{n=0}^N \frac{\lambda^n}{n!}}\right)\\
&>e^{-\lambda}\, \exp \left(\dfrac{ \displaystyle \lambda\  \frac{\lambda^N}{N!}}{\displaystyle \sum_{n=0}^\infty \frac{\lambda^n}{n!}}\right).
\end{align*}
The last expression equals to $e^{-\lambda}\exp \left(  \lambda\ {\frac{\lambda^N}{e^\lambda N!}}\right)$, which is at least $e^{-\lambda}\ \left(1+\lambda\ {\frac{\lambda^N}{e^\lambda \,N!}}\right)$ since $e^x\geq 1+x$ for $x \in \mathbb{R}$. 
\end{proof}
\noindent
By combining these , we get
\begin{align*}
 \exp \left( \dfrac{-\lambda\, \displaystyle \sum_{n=0}^{N-1} \dfrac{\lambda^n}
        {n!} }{\displaystyle \sum_{n=0}^{N} \dfrac{\lambda^n}
        {n!}}\right)\, \sum_{n=0}^{N-1} \dfrac{\lambda^n}
        {n!} -1 &> e^{-\lambda}\ \left(1+\lambda\ {\frac{\lambda^N}{e^\lambda \,N!}}\right)\, e^\lambda\Big(1 -\frac{\lambda^N}{e^\lambda\,N!}\ \frac{N+1}{N+1-\lambda}\Big) -1\\
        & = \frac{\lambda^N}{e^\lambda\,N!} \left(\lambda - \dfrac{N+1}{N+1-\lambda}\right) - \frac{\lambda \,(N+1)}{N+1-\lambda} \left(\frac{\lambda^N}{e^\lambda\,N!}\right)^2\\
        & = \frac{\lambda^N}{(N+1- \lambda)\,e^\lambda\,N!} \left((N+1) \,(\lambda - 1) - \lambda^2 - \dfrac{\lambda^{N+1}\,(N+1)}{e^\lambda\,N!}\right).
        \end{align*}
Let $h(\lambda, N) \coloneqq (N+1) \,(\lambda - 1) - \lambda^2 - \dfrac{\lambda^{N+1}\,(N+1)}{e^\lambda\,N!}$. 
\begin{itemize}
 \item \textbf{\underline{Claim I}:} The function $h(\lambda, N)$ is increasing in $N$.
\begin{proof} Note that
    \begin{align*}
        h(\lambda, N+1) - h(\lambda, N) & = \dfrac{e^{-\lambda}\,\lambda^{N+1}}{(N+1)!}\, \left( (N+1)^2 - \lambda\, (N+2)\right) + \lambda-1.
    \end{align*}
Since $(N+1)^2 - \lambda\, (N+2) = (N+2)\,\left(\frac{(N+1)^2}{N+2} - \lambda\right) = (N+2) \left(N-\lambda + \frac{1}{N+2}\right)>0$, we have $h(\lambda, N+1) - h(\lambda, N)>0$.

\end{proof}

    \item \textbf{\underline{Claim II}:}   $h(\lambda, 5)>0$ for $\frac{3}{2} \leq \lambda \leq 2$.
    \begin{proof}  Letting $N=5$ in $h(\lambda, N)$, we get
     \begin{align*}
        h(\lambda, 5) &= - \lambda^2 + 6\lambda - 6 -\dfrac{1}{20}e^{-\lambda}\,\lambda^6\\
        & = h_1 (\lambda) + h_2(\lambda),
        \end{align*}
        where $h_1 (\lambda) = - \lambda^2 + 6\lambda - 6  $ and $h_2 (\lambda)= -\dfrac{1}{20}e^{-\lambda}\,\lambda^6$. It is easy to see that $h_1(\lambda)$ is increasing and $h_2(\lambda)$ is decreasing on $[\frac{3}{2}, 2]$. Therefore,
   $$h(\lambda, 5) \geq h_1(3/2) + h_2(2) = \frac{3}{4} - \frac{16}{5 e^2}>0.$$
    \end{proof}
\end{itemize}
With these two claims at hand, we conclude that $h(\lambda,N)>0$ on $[\frac{3}{2} , 2] \times [5, \infty)$.

   \end{proof}
\section*{Acknowledgment}  I would like to thank the anonymous referee for careful reading of the manuscript and helpful comments.

\vskip10mm
\noindent Heshan Aravinda \\
Department of Mathematics \& Statistics \\
Sam Houston State University \\
Huntsville, TX 77340, USA \\
heshap@shsu.edu


\begin{thebibliography}{30}
\bibitem{HesMarMel} Aravinda, H., Marsiglietti, A. and Melbourne, J. {\it Concentration inequalities for ultra log-concave distributions.} Studia Mathematica, 265 (2022), 111-120.


\bibitem{BDP} Badiella, L., del Castillo, J. and Puig, P. {\it Ultra log-concavity of discrete order statistics}. Statistics \& Probability Letters, 201, 109900, 2023.

\bibitem{H} Hoeffding, W. {\it On the distribution of the number of successes in independent trials}. The Annals of Mathematical Statistics, 713-721, 1956.
   
\bibitem{JMNS}Jakimiuk, J., Murawski, D., Nayar, P. and Słobodianiuk, S., {\it Log-concavity and discrete degrees of freedom}. Discrete Mathematics, 347(6), 114020, 2024.

\bibitem{JKM} Johnson, O., Kontoyiannis, I. and Madiman, M. {\it Log-concavity, ultra-log-concavity, and a maximum entropy property of discrete compound Poisson measures.} Discrete Applied Mathematics, 161(9), 1232-1250, 2013.

\bibitem{L} Liggett, T. M. {\it Ultra logconcave sequences and negative dependence}. Journal of Combinatorial Theory, Series A, 79(2), 315-325, 1997.

\bibitem{MarMel} Marsiglietti, A. and Melbourne, J.
{\it Geometric and functional inequalities for log-concave probability sequences.} Discrete Comput Geom (2023), https://doi.org/10.1007/s00454-023-00528-7



 \bibitem{MM2} Marsiglietti, A. and Melbourne, J. {\it Moments, Concentration, and Entropy of Log-Concave Distributions.} Preprint, arXiv:2205.08293, 2022.

 \bibitem{Pem}Pemantle, R.
 {\it Towards a theory of negative dependence.} 
J. Math. Phys. 41, no. 3 , 1371 - 1390, 2000.

\bibitem{WW} Xia, W. and Lv, W., {\it Log-concavity and relative log-concave ordering of compound distributions}. Probability in the Engineering and Informational Sciences, 38(3): 579-593. doi:10.1017/S0269964823000293, 2024.
 
\end{thebibliography}
\end{document}